\documentstyle[amssymb,amsfonts]{amsart}

\def\Vol{\hbox{\rm Vol}}
\def\Var{\hbox{\bf Var}}
\def\Conv{\hbox{ Conv}}

\def\P{{\hbox{\bf P}}}
\def\E{{\hbox{\bf E}}}

 at 10 true pt

\def\be#1{ \begin{equation}\label{#1} }

\def\bas{\begin{align*}}
\def\eas{\end{align*}}
\def\bi{\begin{itemize}}
\def\ei{\end{itemize}}

\newenvironment{proof}{\noindent {\bf Proof} }{\endprf\par}
\def \endprf{\hfill  {\vrule height6pt width6pt depth0pt}\medskip}
\def\emph#1{{\it #1}}
\def\textbf#1{{\bf #1}}

\def\BE{{\mathbf E}}

\def\BP{{\mathbf P}}

\def\CK{{\mathcal K}}

\def\BBR {{\mathbb R}}

\def\ep{{\epsilon}}

\def\hs{\hfill $\square$}
\def\dd{\partial}

\theoremstyle{plain}
  \newtheorem{theorem}[subsection]{Theorem}
  \newtheorem{conjecture}[subsection]{Conjecture}

  \newtheorem{proposition}[subsection]{Proposition}
  
  \newtheorem{lemma}[subsection]{Lemma}
  \newtheorem{corollary}[subsection]{Corollary}

\theoremstyle{remark}
  \newtheorem{remark}[subsection]{Remark}

\theoremstyle{definition}

\include{psfig}

\begin{document}

\title{Central limit theorems for   random polytopes in a smooth convex set}

\author{Van Vu}
\address{Department of Mathematics, UCSD, La Jolla, CA 92093-0112}
\email{vanvu@@ucsd.edu}

\thanks{V. Vu is an A. Sloan  Fellow and is supported by an NSF Career Grant.}

\begin{abstract} Let $K$
be a smooth  convex set with volume one in $\BBR^d$. Choose $n$
random points in $K$ independently according to the uniform
distribution. The convex hull of these points, denoted by $K_n$,
is called a {\it random polytope}. We prove that several key
functionals of $K_n$ satisfy the central limit theorem as $n$
tends to infinity.
\end{abstract}

\maketitle

\section {Introduction}

 Let $K$
be a convex set with volume one in $\BBR^d$. Choose $n$ random
points in $K$ independently according to the uniform distribution.
The convex hull of these points, denoted by $K_n$, is called a
{\it random polytope}.
 The study of the key functionals (such as
the volume, the number
 of vertices etc) of $K_n$,
started by Efron \cite{Efr} and R\'enyi and Sulanke \cite{RS},
 is a classical topic in convex geometry (e.g., see  \cite{WW}
 for a survey).

 Let $\Vol (K_n)$ denote the volume of $K_n$. A well known conjecture in the filed
 is
  that this random variable satisfies the
 central limit theorem as $n$ tends to infinity.

\begin{conjecture} \label{conj:CLT} There is a function $\ep(n) $
tending to zero with $n$ such that for every $x$

$$\Big| \P \Big( \frac{ \Vol (K_n) -\E( \Vol (K_n)) }{ \sqrt {\Var (\Vol(K_n))} }
\le x \Big) - \Phi(x) \Big| \le \ep(n)  , $$

\noindent where $\Phi$ denotes the distribution function of the
normal distribution. \end{conjecture}

The conjecture  has been verified only for the case when $K$ is a
ball in $\BBR^2$ \cite{Hsing}.

In the study of  random polytopes in a convex body $K$, the
surface of $K$ plays a critical role. A significant part of the
literature focus on the  following two cases:

\begin{itemize}

\item $K$ is smooth, i.e., the boundary of $K$  is twice
differentiable with positive curvature bounded away from zero and
infinity.

\item $K$ is a polytope.

\end{itemize}

In this paper, we are going to prove the conjecture for the case
when $K$ is a smooth convex set.

\begin{theorem} \label{theo:CLT} Let $K$ be a smooth convex body with volume one in $\BBR^d$. There is a function $\ep(n) $
tending to zero with $n$ such that for every $x$

$$\Big| \P \Big( \frac{ \Vol (K_n) -\E( \Vol (K_n)) }{ \sqrt {\Var (\Vol(K_n))} }
\le x \Big) - \Phi(x) \Big| \le \ep(n)  , $$

\noindent where $\Phi$ denotes the distribution function of the
normal distribution.

\end{theorem}

Using the same method we can  prove similar results for other
functionals of $K_n$. Let $f_i(K_n)$  be the number of
$i$-dimensional facets of $K_n$, where $0 \le i \le d-1$. For
$i=0$,  $f_0(K_n)$ is the number of vertices of $K_n$.

\begin{theorem} \label{theo:CLT3} Let $K$ be a smooth convex body with volume one in $\BBR^d$. There is a function $\ep(n) $
tending to zero with $n$ such that for every $x$ and every $0\le i
\le d-1$

$$\Big| \P \Big( \frac{ f_i (K_n) -\E( f_i (K_n)) }{ \sqrt {\Var (f_i(K_n))} }
\le x \Big) - \Phi(x) \Big| \le \ep(n)  , $$

\noindent where $\Phi$ denotes the distribution function of the
normal distribution.

\end{theorem}

\begin{remark} \label{exponent}  In both theorems one can take $\ep(n)  = n^{-1/(d+1)+o(1)}$. See Remark \ref{exponent1}.    \end{remark}

Our approach is also applicable  for the case when $K$ is a
polytope (pending on the availability of certain results). In
fact, it can also be applied to other models of random polytopes;
see Section 11) for a further discussion.

\vskip2mm

There are  two main technical ingredients in our arguments. The
first is recent results of Reitzner \cite{Rei}  which established
central limit theorems for Poisson point processes. The second is
(also recent) tail estimates by the author \cite{Vu1}.

The paper is organized as follows. In the next section, we
introduce the results from \cite{Rei}. In Sections 3 and 4, we
describe our general approach. In Section 5, we present the above
mentioned tail estimates and a few corollaries. In Section 6, we
prove a tail estimate on the difference $\Vol(K_{n'}) -
\Vol(K_n)$, where $n'$ and $n$ are two large numbers which are
relatively close to each other. As a corollary, we obtain a new
estimate about $\E(\Vol(K_{n'}) - \Vol(K_n))$. The proof of
Theorem \ref{theo:CLT} comes in Section 7. The proof of Theorem
\ref{theo:CLT3} appears in Sections 8 and 9. This latter proof is
somewhat more complex due to a subtle geometrical fact (see the
paragraph following Theorem \ref{theo:Vu2}).  In Section 10, we
summarize recent progresses on the distribution of $\Vol(K_n)$ and
$f_i(K_n)$. In Section 11, we discuss further applications of our
method (for instance, the case when $K$ is a polytope).

In the above sections we sometime omit technical details in order
to maintain the flow of  ideas. These details will be presented in
the appendix at the end of the paper.

\vskip2mm

{\bf \noindent Notations.} In the whole paper, we assume that $n$
is large, whenever needed. The asymptotic notations are used under
the assumption that $n \rightarrow \infty$. The hidden constants
in $O, \Omega$ and $\Theta$  depend on the fixed convex body $K$.
$\Vol$ and $\BP$ denote volume and probability, respectively.
Consider a  (measurable) subset $S$ of $K$ and a random point $x$

$$\Vol (S) = \BP (x \in S). $$

$\BE$ and $\Var$  denote expectation and variance, respectively;
$\log$ denotes the logarithmic with natural base. For a  point set
$P$, $\Conv (P)$ is the convex hull of $P$.

\section {Central limit theorems for Poisson point processes}

Let $X(n)$ be a Poisson point process in $\BBR^d$ of intensity
$n$. Define the random polytope $\Pi_n$ as the convex hull of the
intersection of $K$ with $X(n)$ (a precise  definition  is given
in the next section). In a recent remarkable paper \cite{Rei},
Reitzner proved

\begin{theorem} \label{theo:reitzner1} Let $K$ be a smooth convex body with volume one in $\BBR^d$.
There is a function $\ep(n) $ tending to zero with $n$ such that
for every $x$

$$\Big| \P \Big( \frac{ \Vol (\Pi_n) -\E( \Vol (\Pi_n)) }{ \sqrt {\Var (\Vol(\Pi_n))} }
\le x \Big) - \Phi(x) \Big| \le  \ep(n) . $$

\end{theorem}

\begin{theorem} \label{theo:reitzner2} Let $K$ be a smooth convex body with volume one in $\BBR^d$. Then
 for any $0\le i \le d-1$ there is  a function $\ep_i(n) $ tending to zero with $n$ such that
for every $x$

$$\Big| \P \Big( \frac{ f_i (\Pi_n) -\E( f_i (\Pi_n)) }{ \sqrt {\Var (f_i (\Pi_n))} }
\le x \Big) - \Phi(x) \Big| \le  \ep_i(n)  $$

\end{theorem}

\begin{remark} In Theorems \ref{theo:reitzner1} and \ref{theo:reitzner1}, one can take $$ \ep(n) = c n^{-1/2+ 1/(d+1)} \log
^{2+ 2/(d+1)} n$$ and $$\ep_i(n) = c n^{-1/2+ 1/(d+1)} \log ^{2+
3i+ 2/(d+1)} n, $$  respectively \cite{Rei}. \end{remark}

In the same paper, Reitzner also determined the right order of
magnitude of the variance of $\Vol(K_n)$ and $f_i(K_n)$.

\begin{theorem} \label{theo:reitzner3} Let $K$ be a smooth convex body with volume one in $\BBR^d$. Then
$\Var (\Vol(K_n)) = \Theta (n^{-(d+3)/(d+1)} ) $ and $\Var
(f_i(K_n)) =\Theta (n^{(d-1)/(d+1)}). $
 The same estimates hold for $\Pi_n$.

\end{theorem}

The upper bound for the variance was actually proved in an earlier
paper \cite{Rei2}. A different proof was given by the current
author in \cite{Vu1}.  In \cite{Rei}, Reitzner proved a matching
lower bound.

The proof  in \cite{Vu1} also gives an upper bound for the moment
of any fixed order. A matching lower bound follows from Theorems
\ref{theo:CLT} and \ref{theo:CLT3}.  For details see Section 10.

\section{Passing from the  Poisson model to  the uniform model}

The heart of this paper is  a technique which enables one to pass
a statement for the Poisson model $\Pi_n$ to a similar statement
for the uniform model $K_n$. We will first focus on Theorem
\ref{theo:CLT} concerning the volume.

 In order to prove Theorem \ref{theo:CLT}, given the results in the previous section, it suffices to show
that $\Vol(\Pi_n)$ and $\Vol(K_n)$ have (approximately) the same
expectation and variance. That is exactly what we are going to do

\begin{proposition} \label{pro:expectation} Let $K$ be a smooth convex body with volume one in $\BBR^d$.
Then there is a number  $\delta=1/(d+1) +o(1) $ such that the
following holds.

$$ \E (  \Vol( \overline \Pi_n)) = \E (  \Vol( \overline K_n)) (1+O(n^{-\delta})),
$$

\noindent where $\overline \Pi_n = K \backslash \Pi_n$ and
$\overline K_n = K \backslash K_n$.
\end{proposition}

\begin{proposition} \label{pro:variance} Let $K$ be a smooth convex body with volume one in $\BBR^d$.
 Then there is a number $\delta= 1/(d+1)+o(1) $ such that the
following holds.

$$ \Var (\Vol(\Pi_n)) = \Var (\Vol(K_n)) (1+O(n^{-\delta})). $$
\end{proposition}

\begin{remark} \label{exponent1}
Using the triangle inequality, Theorem \ref{theo:CLT} follows
immediately from these propositions, Theorem \ref{theo:reitzner1}
and the following lemma from \cite{Rei}. Each of these results has
an error term, but the error term $n^{-\delta}$ in the above
proposition turns out to be the dominating one and it leads to
Remark \ref{exponent}
\end{remark}

\begin{lemma} \label{lemma:reitzner}  Let $K$ be a smooth convex body with volume one in
$\BBR^d$. Then

$$\Big| \P(\Vol(K_n) \le x) -\P(\Vol(\Pi_n) \le x) \Big| =O(
n^{-2/(d+1)} \log^{2/(d+1)} n). $$

\end{lemma}

\begin{remark}  Proposition
\ref{pro:expectation} has  already been  verified in \cite{Rei}.
Our proof here is different and works for every  convex body $K$,
without the smoothness assumption. It is mentioned in \cite{Rei}
that it has been conjectured that the variance of $\Vol (\Pi_n)$
and that of $\Vol (K_n)$ are asymptotically the same, but no proof
is known. Proposition \ref{pro:variance} verifies this conjecture
in a strong form.
\end{remark}

In the following, we are going to focus on Proposition
\ref{pro:variance}.  Proposition \ref{pro:expectation} follows
easily from the proof of Proposition \ref{pro:variance}.

Let us now give a better  definition for the Poisson random
polytope $\Pi_n$, mentioned in the previous section. In order to
generate $\Pi_n$, one first generates a random number $n'$ with
respect to the Poisson distribution with mean $n$. Next, one
generates a set $P'$ of $n'$ random points in $K$ with respect to
the uniform distribution. The convex hull of $P'$ is $\Pi_n$.

It is well known from properties of the Poisson distribution that
with high probability $n'$ is close to $n$. To be precise,

$$\P( |n'-n| \ge A \sqrt{n \log n} ) \le n^{-A/4} $$

\noindent for any constant $A$ ($4$ can be replaced by a better
constant but this is not essential). Thus one may think about the
distribution of $\Pi_n$ as a distribution over $K_{n'}$, where $n'
\in [n - A \sqrt{n \log n}, n + A \sqrt{n \log n}]$, for some
large $A$.

\begin{remark} We will prove rigorously that the contribution of those $n'$ outside this interval is
negligible. However, the reader can convince himself (herself)
quickly on this point by observing that every quantity under
investigation  (such as the expectation or variance) is of order
$\Omega(n^{-c})$ for some fixed $c$. Thus by setting $A$
sufficiently large, omitting the $n'$ outside the interval should
not change our analysis significantly. We will use the phrase
"with very high probability" to mean an event which holds with
probability $1-n^{-C}$ where we can set $C$ arbitrarily large.
\end{remark}

Our plan is to show that for any $n'$ in the interval, the
variance of $\Vol (K_{n'})$ and that of $\Vol (K_n)$ is very close
to each other.

\begin{lemma} \label{lemma:main1} There is a number
$\delta=1/(d+1)+o(1)$ such that for any pair $(n,n')$ where $n'= n
+O(\sqrt{n \log n})$ and $n$ is sufficiently large

 $$ \Var (\Vol(K_{n'})) = \Var (\Vol(K_n)) (1+O(n^{-\delta})).
$$ \end{lemma}

This, together with some simple arguments, implies Proposition
\ref{pro:variance}.

\section{Attacking Lemma \ref{lemma:main1}}

We now describe the main idea of how to verify Lemma
\ref{lemma:main1}. This idea is quite general and  works for other
functionals (such as the number of vertices) as well.

  Let us
consider a number   $n'$ in the interval $[n-A \sqrt{n \log n}, A
+ \sqrt{n \log n}]$ for some large constant $A$. Assume (without
loss of generality) that $n'
>n$. Let $\Omega$ denote the product space $K^n$ (equipped with
the natural product measure). A point $P$ in $\Omega$ is an
ordered set $(x_1,...,x_n)$ of $n$ random points (we generate the
points one by one).  The $x_i$ are the coordinates of $P$. We use
$Y(P)$ to denote the volume of the convex hull of $P$ and $\mu$ to
denote the expectation of $Y(P)$.

\begin{remark} $Y(P)$ is, of course,  just another way to express
$\Vol(K_n)$. It is however more convenient to use this notation in
the proof below as it emphasizes the fact that $Y$ is a function
from $\Omega$ to $\BBR$.
\end{remark}

 Define $\Omega', P', \mu'$ similarly (with respect to
$n'$). The variance of $\Vol(K_n)$ is

 $$s= \int_{\Omega} |Y(P)- \mu|^2 \dd P $$

\noindent and  the variance of $\Vol(K_{n'})$ is

 $$s'= \int_{\Omega'} |Y(P')- \mu'|^2 \dd P'. $$

Our goal is to show that $\delta$

$$ |s'-s| \le n^{-\delta} s, $$

\noindent for some $\delta$ as claimed.  We are going to use a
coupling argument. Consider a point $P' =(x_1, \dots, x_{n'})$ in
$\Omega'$ and  the canonical decomposition

$$P' = P \cup Q $$

\noindent where $P=(x_1, \dots, x_n)$ and $Q=(x_{n+1}, \dots,
x_{n'})$.  In order to compare $s$ and $s'$, we rewrite $s$ as

$$s = \int_{\Omega'}  (Y(P) -\mu)^2  \dd
P' $$

\noindent where
 $Y(P)$ is understood as a  function of $P'$ which
depends on the first $n$ coordinates of $P'$.
 We now can write $s'-s$ as

\begin{equation} \label{equa:ss'1} s' -s = \int_{\Omega'}
(Y(P') -\mu')^2 - (Y(P)-\mu)^2 ) \dd P'. \end{equation}

\noindent   Next, we observe that

$$(Y(P')  -\mu')^2 - (Y(P)-\mu)^2 =  \big((Y(P') -\mu') + (Y(P)-\mu) \big) \big((Y(P') -\mu') - (Y(P)
-\mu ) \big). $$

Using a recent tail estimate  in \cite{Vu1}, we can   prove that
with very high probability, the quantities $|Y(P') -\mu'|$ and
$|Y(P)-\mu|$ are small, of order $O(\sqrt s \log^{O(1)} n)$ in
particular. This implies that typically

\begin{equation} \label{equ:sum}   Y(P') -\mu') + (Y(P)-\mu)  = O(\sqrt s \log^{O(1)} n)
\end{equation}

\noindent and

\begin{equation} \label{equ:difference}  (Y(P')  -\mu') - (Y(P) -\mu )  =  O(\sqrt s
\log^{O(1)} n).
\end{equation}

\noindent These two inequalities, however, are not enough for  our
purpose. They only imply that the product of the left hand sides
is of order $O(s \log^{O(1)} n)$, while we want something
significantly small than $s$ (the exponent in the logarithmic term
is, unfortunately, positive). While the estimate for the sum $
(Y(P') -\mu') + (Y(P)-\mu) $ seems to be near optimal, we can
prove that  a much better estimate holds for the difference $
(Y(P') -\mu') - (Y(P) -\mu )$. In order to achieve this, we
rewrite this difference as follows

$$ (Y(P') - Y(P)) -(\mu'- \mu). $$

\noindent We are going to show that with very high probability

\begin{equation} \label{equ:differencenew} Y(P') - Y(P) = O( \sqrt s
n^{-\delta}).
\end{equation}

\noindent  This will implies that $\mu'- \mu$ is of the same
order. Thus, we gain a significant extra term $n^{-\delta}$ as
claimed. In this part of the proof, the facts that $P'$ contains
$P$  and $n'-n =O(\sqrt{n \log n})$ will play essential roles.

\section{Concentration of the Volume}

Again we use short hand $Y$ to denote the random variable
$\Vol(K_n)$.  Recall that the variance of $Y$ is $\Theta
(n^{-(d+3)/(d+1)})$. One expects that $Y$ has sub-Gaussian tail,
namely

$$ \P(|Y - \E(Y)| \ge \sqrt{\lambda n^{-(d+3)/(d+1)} }) \le
\exp(-c\lambda) $$

\noindent for some constant $c$ and a large range of $\lambda$.
This intuition  was confirmed in \cite{Vu1}. The following theorem
is equivalent to  Corollary 2.8 from that paper.

\vskip2mm

\begin{theorem} \label{theo:Vu1}  Let $K$ be a smooth convex body with volume one in
$\BBR^d$. There are positive constants $c$ and $ \alpha$ such that
the following holds. For any $0 < \lambda \le  \alpha
n^{\frac{(d-1)(d+3)}{(d+1)(3d+5)}}$

\begin{equation} \label{equa:sharpconcentration} \BP(|Y-\BE(Y)| \ge \sqrt {\lambda n^{-(d+3)/(d+1)} }) \le
2\exp(-c\lambda) + \exp (-c n^{(d-1)/(3d+5)}). \end{equation}

\end{theorem}

It is useful to compare this theorem with the central limit
theorem. In the central limit theorem, we  have, for sufficiently
large $\lambda$

\begin{equation} \label{equa:CLT} \BP(|Y-\BE(Y)| \ge \sqrt {\lambda \Var(Y)}
)     \approx \frac{2}{\lambda} \exp(-\lambda/2). \end{equation}

Since $\Var(Y) =\Theta(n^{-(d+3)/(d+1)})$, the two bounds are
comparable. The key difference here is that in the central limit
theorem,  $\lambda$ has to be fixed and one takes $n$ to infinity,
while \eqref{equa:sharpconcentration} also  works for  $\lambda$
depending on $n$.

\begin{remark} For the purpose of proving Theorem \ref{theo:CLT}, it  actually
suffices to use a weaker theorem (Theorem 2.1 of \cite{Vu1}) whose
proof is significantly simpler. In this weaker theorem  instead of
$\sqrt{\lambda n^{-(d+3)/(d+1)} }$ one has $\sqrt{\lambda
n^{-(d+3)/(d+1)} \log n }$ (an extra logarithmic term). We state
the above more precise theorem because we will discuss its
analogue for $f_i(K_n)$ later. \end{remark}

We are also going to use the following standard fact, which can be
routinely verified using the Laplace transform or properties of
binomial coefficient (see \cite{JLR}, Chapter 2).

\begin{lemma} \label{lemma:chernoff} Let $X$ be the sum of $m$
i.i.d. indicator random variables. Then for any $t \ge 0$

$$\P (X \ge \E(X)+ t) \le \exp(- \frac{t^2}{ 2(\E(X) + t/3) } ). $$

\end{lemma}

\section{Differences of Volumes}

Consider two integers $n > n'$. We generate a set $P'$ of $n'$
random points as follows. First we generate a set $P$ of $n$
random points and next generate a set $Q$ of $m=n-n'$ random
points and then define $P'= P \cup Q$. Our goal is to show that if
$n'-n$ is small compared to $n$, then with very high probability,
$Y(P') -Y(P)$ is small. The argument in this section works for
every convex body $K$, {\it without}  the smoothness assumption.

 For a convex body $K$ and a half space $H$,  we
 call the intersection $H \cap K$ an {\it $\ep$-cap} if $\Vol
(H \cap K) =\ep$. The union of all $\ep$-caps is the {\it
$\ep$-wet part} of $K$. The complement of the $\ep$-wet part is
called the {\it $\ep$-floating body} of $K$. We denote the
$\ep$-floating body and the $\ep$-wet part by $F_{\ep}$ and
$\overline F_{\ep}$, respectively. The volume of the $\ep$-wet
part play an important role and we denote it by  $\rho_{\ep}$.

\vskip2mm

 Consider two points $x$ and $y$ in
$\overline {F_{\ep}}  $. We say that $x$ {\it sees} $y$ (with
respect to $F_{\ep}$)  if the segment $xy$ does not intersect the
floating body $F_{\ep}$. Let $S_{x, \ep}$ denote the set of those
$y$ that $x$ sees. Define

$$g(\ep) = \sup_{x \in\overline {F_{\ep}}   } \Vol (S_{x, \ep}). $$

\noindent The set $S_{x, \ep}$ is  the union of all $\ep$-caps
containing $x$. The quantity which is important to us is

$$\ep^{\ast} :=  \nu \log n /n $$

\noindent where $\nu$ is a sufficiently large constant.

\begin{lemma} \label{lemma:difference} Let $K$ be a convex body of volume one in $\BBR^d$.
For any constant $\beta$, there are constants $\gamma$ and $\nu$
such that for $\ep^{\ast} =\nu \frac{\log n}{n}$, the following
holds

$$\P \Big(Y(P') -Y(P) \ge  (\gamma+1) m \rho (\ep^{\ast}) g(\ep^{\ast})
+ \gamma g(\ep^{\ast} ) \log n  \Big)  \le n^{-\beta}. $$

\end{lemma}

\begin{corollary}\label{cor:difference} Let $K$ be a convex body
of volume one in $\BBR^d$ and $\mu(n)$ the expectation of the
volume of $K_n$. Then

\begin{equation} \label{cor:difference1} \mu (n+m) - \mu (n) = O(
g(\ep^{\ast}) m \rho (\ep^{\ast}) + g(\ep^{\ast}) \log n ).
\end{equation}

\noindent If $K$ is smooth, then

\begin{equation} \label{cor:difference2} \mu(m+n) - \mu(n) = O( (m n^{-2/(d+1)}+1) n^{-1} \log^{O(1)} n).
\end{equation}

\noindent If $K$ is a polytope, then

\begin{equation} \label{cor:difference3} \mu(m+n) - \mu(n) = O( (m n^{-1}+1)n^{-1}  \log^{O(1)} n). \end{equation}

\end{corollary}

\begin{remark} Using a stronger definition of smoothness, one
 can slightly improve \eqref{cor:difference2}. Indeed,
if one assumes that $K \in \CK^{+}_k$, then it is known that

$$\mu_n =\sum_{i=1}^k c_i n^{-(i+1)/(d+1)} + o( n^{-(k+1)/(d+1)})
$$

\noindent where $c_i$ are constants depending on $K$. Thus by
taking $k$ sufficiently large and using Taylor's expansion, one
can prove \eqref{cor:difference2} without the logarithmic term.
Here $\CK^{+}_k$ consists of convex body whose boundary is $k$
time differentiable.

\end{remark}

\begin{proof} In order to prove Lemma \ref{lemma:difference}, we
are going to need the following lemma from \cite{BD}. (For a
different proof of this lemma using the VC dimension,  see Section
4 of \cite{Vu1}.)

\begin{lemma} \label{lemma:floating} There are  positive constants $c$ and $c'$
 such that the following
holds for every sufficiently large $n$. For any $\ep \ge c' \log
n/n$, the probability that $K_n$ does not contain $ F_{\ep}$ is at
most $\exp (-c \ep n ).$

\end{lemma}

Let $E$ be the event that $F_{\ep^{\ast} }$ is not contained in
$\Conv (P)$.  For any $T$

$$ \P (Y(P') -Y(P) \ge T)  \le \P(E)+ \P( Y(P') -Y(P) \ge T | \overline E).
$$

By setting $\nu$ in the definition of $\ep^{\ast}$ we can assume
that $\ep^{\ast}$ satisfies the condition of Lemma
\ref{lemma:floating} and the probability $\exp(-c \ep^{\ast}  n)$
we get from this lemma is at most $n^{-\beta-1}$. Thus we have

\begin{equation} \label{equa:beta1}  \P (Y(P') -Y(P) \ge T)  \le n^{-\beta-1} + \P( Y(P') -Y(P) \ge T | \overline
E). \end{equation}

Observe that when $\overline E$ holds (namely, $F_{\ep^{\ast}}
\subset \Conv (P)$), the volume of $\Conv(P)$ increases by at most
$g(\ep^{\ast})$ when we add to $P$ a new point from $K \backslash
\Conv(P)$. Since $F_{\ep^{\ast}} \subset \Conv (P)$, $K \backslash
\Conv(P) \subset \overline F_{\ep^{\ast}}$ and then

$$ \P( Y(P') -Y(P) \ge T | \overline E) \le \P(|Q \cap \overline
F_{\ep^{\ast}}| \ge T/g(\ep^{\ast})) . $$

\noindent By substituting

$$T= g(\ep^{\ast}) \Big( (\gamma+1)  m \rho_{\ep^{\ast}}  +  \gamma \log n \Big), $$ we have

$$ \P(|Q \cap \overline F_{\ep^{\ast}}| \ge T/g(\ep^{\ast})) = \P \Big(|Q \cap \overline
F_{\ep^{\ast}}| \ge  (\gamma  +1)m \rho_{\ep^{\ast}}  +  \gamma
\log n) \Big).
$$

\noindent The quantity $|Q \cap \overline F_{\ep^{\ast}}|$ is a
sum of $m$ i.i.d. indicator random variables (each of which
indicates the event that a random point belongs to the wet-part
$F_{\ep^{\ast}}$ or not) with expectation

$$\Vol (\overline F_{\ep^{\ast}} = \rho_{\ep^{\ast}}  . $$

\noindent Thus, the expectation of $|Q \cap \overline
F_{\ep^{\ast}}|$ is $m \rho_{\ep^{\ast}} $. By applying Lemma
\ref{lemma:chernoff} with

$$t= \gamma (m \rho_{\ep^{\ast}}  +\log n), $$

\noindent  we conclude that the probability in question is bounded
from above by

$$\exp\Big(- \frac{ \gamma^2 (m \rho_{\ep^{\ast}}  +\log n)^2} {2( m
\rho_{\ep^{\ast}}   + \frac{1}{3} \gamma (m \rho_{\ep^{\ast}}
+\log n))} \Big). $$

\noindent It is simple to prove that for any $\gamma \ge 3$, the
exponent

$$ \frac{ \gamma^2 (m \rho_{\ep^{\ast}}  +\log n)^2} {2( m
\rho_{\ep^{\ast}}   + \frac{1}{3} \gamma (m \rho_{\ep^{\ast}}
+\log n))} $$

\noindent is bounded from above (with room to spare)  by
$\gamma/2$. Set $\gamma $ large enough so that $\gamma/2 > \beta
+1$. We have

\begin{equation} \label{equa:beta2} \exp\Big(- \frac{ \gamma^2 (m \rho_{\ep^{\ast}}  +\log n)^2}
{2( m \rho_{\ep^{\ast}}   + \frac{1}{3} \gamma (m
\rho_{\ep^{\ast}} +\log n))} \Big ) \le n^{-\gamma/2}  \le
n^{-\beta -1} .
\end{equation}

\noindent Since $2n^{-\beta-1} < n^{-\beta}$, \eqref{equa:beta1}
and  \eqref{equa:beta2} conclude the proof. \end{proof}

\section {Proof of Proposition \ref{pro:variance}}

We first prove that there is a constant $\delta >0$ such that for
any $n'$ where $|n'-n| = O(\sqrt {n \log n})$,

$$ \Var (\Vol(K_{n'})) = \Var (\Vol(K_n)) (1+O(n^{-\delta})). $$

Without loss of generality, we can assume that $n'>n$. Let $C$ be
a large constant. We call a point $P' =(x_1,..., x_{n'}) \in
\Omega'$ $C$-{\it typical} if the following hold

\begin{itemize}

\item
$|Y(P') - \mu'| \le  C  \sqrt {n^{-(d+3)/(d+1)} \log n} $

\item
$   |Y(P) - \mu| \le  C  \sqrt {n^{-(d+3)/(d+1)} \log n} $

\item $Y(P') -Y(P) \le Cg(\ep^{\ast}) ( m \rho_{\ep^{\ast}}  +
\log n ) $

\end{itemize}

The results in the previous two sections show that there is a
positive  constant $\alpha$ independent of $C$ such that  $P'$ is
$C$-typical with probability at least $1- n^{-\alpha C}$. Since
$K$ is smooth, $g(\ep)= O (\ep)$ and $\rho(\ep) = \Theta
(\ep^{2/(d+1)})$  for all sufficiently small $\ep$. Substituting
$\ep^{\ast} = \Theta (\log n/ n)$ and  $m = O(\sqrt{n \log n} ) $,
we obtain

$$ Cg(\ep^{\ast}) ( m \rho_{\ep^{\ast}}  + \log n )  = O(
n^{-(d+5)/2(d+1)} \log^{O(1)} n ). $$

\noindent This implies that  if  $P'$ is $C$-typical then

\begin{align*} |(Y(P')  -\mu')^2 - (Y(P)-\mu)^2 | &= O \Big(   \sqrt
{ n^{-(d+3)/(d+1) }\log n } \times n^{-(d+5)/2(d+1)} \log^{O(1)} n
\Big)
\\ &=O \Big(  n^{-(d+4)/(d+1)} \log ^{O(1)} n \Big). \end{align*}

\noindent   On the other hand, the variance $s$ is of order $
\Theta (n^{-(d+3)/(d+1)}). $ Thus, we obtain the claim of Lemma
\ref{lemma:main1} with $\delta= 1/(d+1) -o(1)$.

\vskip2mm

Now we are going to conclude the proof of Proposition
\ref{pro:variance}. Consider the interval $I= [n - A\sqrt{n \log
n}, n + A \sqrt{n \log n}]$. For every number ${n'} $ in this
interval, let $E_{n'} $ denote the events that $n'$ is sampled
(according to the Poisson distribution with mean $n$) and $E_0$
denote the event that  the sampled number does not belong to the
interval. We have

$$\Var (\Pi_n) = \E _{n'}  (\Var (\Vol(\Pi_n)|E_{n'} )) + \Var (\E(\Vol
(\Pi_n|E_{n'} )), $$

\noindent where $n' \in I$ or $n'=0$. By increasing $A$, we can
assume that $\P(E_0) \le n^{-B}$ for a sufficiently large $B$.
Thus, the contribution of $\Vol(\Pi_n|E_0)$ is negligible.
Moreover, for any ${n'} \in I$,

$$\Vol(\Pi_n)|E_{n'}  =\Vol (K_{n'} ). $$

\noindent Thus, we have

$$ \Var (\Vol(\Pi_n)|E_{n'} ) = \Var (\Vol (K_{n'})) =
\Var(\Vol(K_n)) (1+ O(n^{-\delta})) $$

\noindent     and also

$$ \E(\Vol (\Pi_n|E_{n'} )) = \E(\Vol(K_{n'})) =\E(\Vol(K_n)) (1+
n^{-\delta}). $$

\noindent These imply that

$$\Var(\Vol (\Pi_n)) = \Var (\Vol (K_n)) (1+ O(n^{-\delta})), $$

\noindent proving Proposition \ref{pro:variance}. \hs

\section {Proof of Theorem \ref{theo:CLT3}}

We first prove this theorem for the  case $i=0$ ($f_0(K_n)$ counts
the number of vertices of $K_n$).  The (slightly more technical)
proof for a  general  $i$ follows in the next section. We need to
prove the following proposition.

\begin{proposition} \label{pro:variance2} Let $K$ be a smooth convex body with volume one in $\BBR^d$.
 Then there is a number  $\delta=1/(d+1) +o(1) $ such that the
following holds.

$$ \Var (f_0(\Pi_n)) = \Var (f_0(K_n)) (1+O(n^{-\delta})). $$
\end{proposition}

Once the proposition is available, one  can use  Theorem
\ref{theo:reitzner2} to conclude the proof. The above  proposition
follows from a variant of Lemma \ref{lemma:main1}.

\begin{lemma} \label{lemma:main2} There is a number
$\delta =1/(d+1)+o(1)$ such that for any pair $(n,n')$ where $n'=
n +O(\sqrt{n \log n})$ and $n$ is sufficiently large

 $$ \Var (f_0(K_{n'})) = \Var (f_0(K_n)) (1+O(n^{-\delta})).
$$ \end{lemma}

In order to prove this lemma, we use the ideas described in
Section 4.  We are going to use short hand  $Z$ to denote
$f_0(K_n)$. Theorem \ref{theo:Vu1} is replaced by the following
result, which is a corollary of Theorem 2.11 from \cite{Vu1}.

\begin{theorem} \label{theo:Vu2} Let $K$ be a smooth convex body with volume one in
$\BBR^d$. There are positive constants $c$ and $ \alpha$ such that
the following holds. For any $0 < \lambda < \alpha
n^{\frac{d-1}{3(d+1)}}$, we have

$$\BP(|Z-\BE(Z)| \ge \sqrt {\lambda n^{(d-1)/(d+1)} }) \le 2\exp(-c\lambda) +
\exp (-c n^{(d-1)/(3d+5)}). $$

\end{theorem}

While it seems that we have all necessary tools in hand, the proof
of Lemma \ref{lemma:main2} is still not a straightforward
modification of the proof of Lemma \ref{lemma:main1}. In fact this
proof is much more delicate  due to the following subtle
obstruction: $Z$ is not a {\it monotone} random variable. Recall
that when we bound the difference between the volumes  $ Y (P')$
and $ Y (P)$,
 we relied on the following observation: Adding a point
from $P'\backslash P$ to $P$  cannot increase $Y(P)$ by too much,
unless $P$ itself is very distorted, which happens with a
negligible probability. When we repeat this argument for the
number of vertices instead of the volume, we have to take into
account the fact that the number of vertices of $\Conv (P)$ can
actually {\it decrease} when we add a new point to $P$. Therefore,
we have to give a upper bound for both $Z(P') -Z(P)$ and
$Z(P)-Z(P')$.

The handling of $Z(P')-Z(P)$ is similar to that of $Y(P')-Y(P)$.
If $\Conv (P)$ contains the floating body $F_{\ep^{\ast}}$, then a
point from $Q= P' \backslash P$ would increase $Z(P)$ by at most
one if it falls in the wet part $\overline F_{\ep^{\ast}}$. Thus,
$Z(P')-Z(P)$ (conditioned on $P$ contains $F_{\ep^{\ast}}$) is
bounded from above by the sum of $m=|Q|$ i.i.d. indicator random
variables with expectation $\rho_{\ep^{\ast}} = \Vol(\overline
F_{\ep^{\ast}})$. Using Lemma \ref{lemma:chernoff}, it is easy to
show that this sum is

$$O( m \rho_{\ep^{\ast}} + \log n) $$

\noindent with probability $1-n^{-C}$ where one can have the
constant $C$ arbitrarily large.

\vskip2mm

The handling of $Z(P)-Z(P')$ is more delicate, as  a point from
$Q= P' \backslash P$ can decrease  $Z(P)$ by an arbitrary amount.
For $k \ge 2$, we say that a point $x$ is $k$-{\it wide} with
respect to $P$ if $x$ lies outside the convex hull of $P$ and sees
$k$ vertices of this convex hull.

\noindent Define

  $$U_{k, P} =\{x| x \,\,\, \hbox{is} \,\,\,
  k-\hbox{wide with respect to } \,\,\, P \}. $$

The key new ingredient is the following lemma  which asserts that
with high probability, $U_{k,P}$ has small volume.

\begin{lemma} \label{lemma:main3} There are positive constants $c_1,
c_2, c_3, c_4$  such that the following holds. For any $k \ge c_3$
and $T\ge c_3  \max \{ \frac{k}{c_4n}, \rho_{ \frac{k}{c_4n}}
\exp(-c_2 k) \}$ we have

$$ \BP( \Vol  ( U_{k, P} ) \ge T) \le \exp(-c_1nT) $$

\noindent where $P$ is a set of $n$ random points. Furthermore,
for any constant $C$ there is a constant $c_5$ such that with
probability at least $1-n^{-C}$, there is no point in $K$ which is
$c_5 \log n$-wide with respect to $P$.

\end{lemma}

The first half of the lemma is Lemma 7.1 from \cite{Vu1}. In order
to maintain the flow of the presentation, we defer the proof of
the second half to the appendix.

Let $c_6$ be a large positive constant. For each $c_3 \le k \le
c_5 \log n$, define

$$ T_k =\max  \{ c_3  \max \{ \frac{k}{c_4n}, \rho_{\frac{k}{c_4n}} \exp(-c_2 k) \}, c_6 \log n/ n \} .$$

The first statement of Lemma \ref{lemma:main3} implies that with
probability at least

$$1- \sum_{k= c_3}^{c_5 \log n} \exp(-c_1 c_6 \log n )$$ we have

$$ \Vol ( U_{k, P} ) \le T_k $$

\noindent for all $c_3 \le k \le c_5 \log n$.

 We call the random set
$P$ {\it perfect} if the following three conditions hold

\begin{itemize}

\item $\Vol ( U_{k, P} ) \le T_k$ for all $c_3 \le k \le c_5 \log
n$.

\item  $\Vol ( U_{k, P} )=0$ for any $k
> c_5 \log n$.

\item $ \Vol (K \backslash P) \le c_7 \rho_{1/n}.$

\end{itemize}

 The above argument and the second part of Lemma
\ref{lemma:main3} and Theorem \ref{theo:Vu1} together imply that
for any given constant  $C$, we can choose the constants $c_3,
c_5$ and $c_7$ such that $P$ is perfect with probability at least
$1-n^{-C}$, (The values of these constants do not play any role in
what follows.)

Now consider a perfect set $P$. We are going to expose a set $Q$
of $m$ random points and see how it changes the number of vertices
of $\Conv(P)$. Let $x$ be a random point and $S$ be the number of
vertices of $\Conv(P)$ seen by $x$. The distribution of $S$, due
to the fact that $P$ is perfect, satisfies

$$\P( S \ge k)  \le T_k $$

\noindent for $c_3 \le k \le c_5 \log n$ and

$$\P( S \ge k) =0$$

\noindent for $k > c_5 \log n$. Furthermore, for any $1 \le k <
c_3$,

$$\P(S =k) \le \Vol (K \backslash P) \le c_7 \rho_{1/n}. $$

\noindent The difference $X= Z(P)- Z(P')$, conditioned on a
perfect set $P$, is bounded from above by
 the sum of $m$ independent copies of $S$. A routine Laplace
 transform argument gives

 \begin{equation} \label{equ:laplace} \P(X \ge D) = \P(e_0^X \ge e_0 ^D) \le \E(e_0^X) e_0^{-D}  = \E(e_0^{S})^m
 e_0^{-D} \end{equation}

\noindent for any number  $e_0 >1$. It remains to bound
$\E(e_0^{S})$, for some special choice of $e_0$. Observe that

\begin{equation} \label{equa:expectation1} \E(e_0^S) \le 1 + \sum_{k=1}^{\infty} e_0^k \P(S=k) .
\end{equation}

\noindent The information about the distribution of $S$ enables us
to bound this from above by

\begin{equation} \label{equa:expectation2}  1 + O(\rho_{1/n})  + \sum_{k=c_3}^{c_5 \log n} e_0^k T_k.
\end{equation}

\noindent By the definition of $T_k$, there is a constant $c_8$
such that

\begin{equation} \label{equa:expectation3}  \sum_{k=c_3}^{c_5 \log n} e_0^k T_k \le c_8
\sum_{k=c_3}^{c_5 \log n} e_0^k \Big( \frac{k}{c_4n}+
\rho_{\frac{k}{c_4n}} \exp(-c_2 k) +  \frac{\log n}{ n} \Big).
\end{equation}

\noindent Set $e_0$ sufficiently close to one, we have

\begin{equation} \label{e0small} \sum_{k=c_3}^{c_5 \log n} e_0^k \frac{k}{c_4n} \le n^{-.9}
\,\,\hbox{and} \,\, \sum_{k=c_3}^{c_5 \log n} e_0^k \frac{\log
n}{n} < n^{-.9}. \end{equation} Again by setting $e_0$ close to
one, we can assume that $e_0 ^k \exp(-c_2 k) \le \exp(-c_2k/2)$.
Thus,

$$ \sum_{k=c_3}^{c_5 \log n} e_0^k T_k =O\Big( \sum_{k=c_3}^{c_5 \log n} \rho_{ \frac{k}{c_4n}} \exp(-c_2 k/2) + n^{-.9}
\Big). $$

 \noindent Using the facts that $$\rho_{\ep} =\Theta
(\ep^{2/(d+1)}) \,\, \hbox{and} \,\,  \int_{c_3}^{\infty}
x^{2/(d+1)} \exp(-c_2x/2) \dd x \,\, \hbox{ converges},$$  one can
easily see that

\begin{equation} \label{equa:e0small1}  \sum_{k=c_3}^{c_5 \log n} \rho_{ \frac{k}{c_4n}} \exp(-c_2
k/2)=O(n^{-2/(d+1)}) =O(\rho_{1/n}). \end{equation}

\noindent The error term $n^{-.9}$ is negligible. Thus we can
conclude that there is a positive constant $\gamma $ such that for
any $1< e_0 < 1+\gamma$

$$\E(e_0^{S}) = 1+ O(\rho_{1/n}). $$

\noindent By \eqref{equ:laplace}, the probability $\P(X \ge D)$ is
bounded from above by

\begin{equation} \label{equa:e0small2}  \E(e_0^{S})^m e_0^{-D} = \exp( O(m  \rho_{1/n}) -  D \log
e_0) \le \exp(c_9 m \rho_{1/n} -D \log e_0) \end{equation}

\noindent for some constant $c_9$ which depends on the choice of
$\ep_0$. Setting

$$D = \frac{1}{\log e_0} (c_9 m \rho_{1/n} + C \log n)  $$

\noindent (where $C$ is a large constant) we conclude that

$$\P(X \ge D) \le n^{-C}. $$

\noindent The rest of the proof is the same as in the proof for
volume. The value of $D$ is negligible compared to $\sigma$, the
standard deviation of $Z$. Indeed, recall that the latter is
$\Omega (n^{(d-1)/2(d+1)})$ (see \cite{Rei}) and  $m= O(\sqrt{n
\log n})$. Thus,

$$D= O( \sigma n^{-1/(d+1) } \sqrt{\log n}) .$$

\noindent This, in particular, shows that we can set
$\delta=1/(d+1)+o(1)$ as claimed in the theorem. The reader is
invited to fill in the details. \hs

\section{Proof of Theorem \ref{theo:CLT3} (general case)}

In order to obtain Theorem \ref{theo:CLT3} in full generality, we
need to make a few technical modifications. Let us denote by $Z_i$
the number of facets of dimension $i$ of $K_n$ (thus $Z= Z_0$).
Similarly, we define $S_i$ and $X_i$ as the analogues of $S$ and
$X$, respectively. The main reason for modification is that $S_i$
can be significantly larger than $S_0=S$.

To start, we have the following generalization of Theorem
\ref{theo:Vu3}, which we can prove by slightly modifying the proof
of Theorem \ref{theo:Vu3} from \cite{Vu1}.

\begin{theorem} \label{theo:Vu3} Let $K$ be a smooth convex body with volume one in
$\BBR^d$. For any $0 \le i \le d-1$, there are positive constants
$c_i$ and $ \gamma_i$ and $\beta_i$  such that the following
holds. For any $0 < \lambda < \gamma_i n^{\beta_i}$,

$$\BP(|Z_i-\BE(Z_i)| \ge \sqrt {\lambda n^{(d-1)/(d+1)} }) \le 2\exp(-c_i\lambda) . $$

\end{theorem}

We are going to use the following lemma, which is a corollary of
the classical Upper Bound theorem.

\begin{lemma} \label{lemma:upperbound} For every $0 \le i < d$, let $\alpha_i= \min \{i, d-i \}$. Then $$|S_i| =O( |S_0|
^{\alpha_i}). $$ \end{lemma}

\noindent We are going to consider the following  analogue of
\eqref{equ:laplace}

\begin{equation} \label{equ:laplacei} \P(X_i \ge D_i) = \P(e_i^{X_i} \ge
e_i ^{D_i}) \le \E(e_i^{X_i} ) e_i^{-D_i}  = \E(e_i^{S_i})^m
 e_i^{-D_i} \end{equation}

\noindent where $e_i$ and $D_i$ are going to be defined. Similar
to \eqref{equa:expectation1}, we have

\begin{equation} \label{equa:expectation1i} \E(e_i^{S_i} ) \le 1 + \sum_{k=1}^{\infty} e_i^k \P(S_i=k) .
\end{equation}

\noindent Using the information about  the distribution of $S_0
=S$ and Lemma \ref{lemma:upperbound},  we can bound the right hand
side from above by

\begin{equation} \label{equa:expectation2i}  1 + O(\rho_{1/n})  + O( \sum_{k=c_3}^{c_5 \log  n}
e_i^{\beta k^{\alpha_i}}  T_k k^{\alpha_i-1})
\end{equation}

\noindent for some constants $\beta$, $c_3$ and $c_5$.  The key
difference (compared to the case $i=0$) here is that the exponent
of $e_i$ is now a polynomial in $k$.  Thus, we need to find some
suitable values for $e_i$ and $D_i$. We set

$$e_i = 1 + \frac{\ep_i}{ \log^{\alpha_i} n } $$

\noindent for some sufficiently small positive constant $\ep_i$.
 By the definition of $T_k$, there is a constant $c_8$ such that

\begin{equation} \label{equa:expectation3i}  \sum_{k=c_3}^{c_5 \log  n} e_i^{\beta k^{\alpha_i}} T_k k^{\alpha_i-1}  \le c_8
\sum_{k=c_3}^{c_5 \log n} e_i^{\beta k^{\alpha_i}} \Big(
\frac{k}{c_4n}+ \rho_{\frac{k}{c_4n}} \exp(-c_2 k) +  \frac{\log
n}{ n} \Big) k^{\alpha_i-1}.
\end{equation}

\noindent Set the constant $\ep_i$ in the definition of $e_i=
\frac{\ep_i}{\log^{\alpha_i} n}$ sufficiently small, we have

\begin{equation} \label{e0small} \sum_{k=c_3}^{c_5 \log  n} e_i^{\beta k^{\alpha_i} } \frac{k^{\alpha_i} }{c_4n} \le n^{-.9}
\,\,\hbox{and} \,\, \sum_{k=c_3}^{c_5 \log^{\alpha_i}  n}
e_i^{\beta k^{\alpha_i}} k^{\alpha_i-1}  \frac{\log n}{n} <
n^{-.9}.
\end{equation} and also $e_i ^{\beta k^{\alpha_i}} \exp(-c_2 k) \le \exp(-c_2k/2)$.
Thus,

$$ \sum_{k=c_3}^{c_5  n} e_i^{\beta k^{\alpha_i}}  T_k  k^{\alpha_i-1} =
O\Big( \sum_{k=c_3}^{c_5 \log n} \rho_{ \frac{k }{c_4n}} \exp(-c_2
k/2) k^{\alpha_i-1}  + n^{-.9} \Big). $$

\noindent   Using the facts that $$\rho_{\ep} =\Theta
(\ep^{2/(d+1)}) \,\, \hbox{and} \,\, \int_{c_3}^{\infty}
x^{\alpha_i-1+ 2/(d+1)} \exp(-c_2x/2) \dd x \,\, \hbox{
converges},$$ we obtain the following analogue of
\eqref{equa:e0small1}

\begin{equation} \label{equa:eismall1}  \sum_{k=c_3}^{c_5 \log n} \rho_{ \frac{k }{c_4n}} \exp(-c_2
k/2)  k^{\alpha_i-1} =O(n^{-2/(d+1)}) =O(\rho_{1/n}).
\end{equation}

\noindent Thus we have the following analogue of
\eqref{equa:e0small2}

\begin{equation} \label{equa:eismall2}  \P(X_i \ge D_i) \le  \E(e_i^{S_i})^m e_i^{-D_i}
 e_i^{-D_i}  =\exp(- c_9 m \rho_{1/n} - D_i \log e_i). \end{equation}

 \noindent In order to make the right hand side at most $n^{-C}$,
 for a sufficiently large constant $C$,
 we need to set

 $$D = \frac{1}{\log e_i} (c_9  m \rho_{1/n} +  C\log n). $$

 \noindent Recall that $e_i=1+ \Omega (\frac{1}{\log^{\alpha_i}
 n})
 $. Thus,

 $$D= O (m \rho_{1/n} +  \log n) \log^{\alpha_i} n = O(\sigma_i
 n^{-1/(d+1)} \log^{1/2 +\alpha_i} n), $$

 \noindent where $\sigma_i= \Theta (n^{(d-1)/2(d+1)})$ is the
 standard deviation of $Z_i$. This concludes the proof. \hs

 \section{A summary }

In this section we give  a brief summary about  the volume and the
number of facets of a random polytope in a smooth convex body.
Most of these results (except those concerning expectation) have
been obtained in the last three years or so, using tools from
modern probability.
 In this theorem, we use the notation

$$M_k=  \int |\Vol (K_n) - \E(\Vol (K_n)| ^k .$$

\noindent For $k$ even, $M_k$ is the $k$th moment.

 \noindent {\it Volume.}

 \begin{theorem} \label{theo:sumarizevolume} Let $K$ be a smooth convex body of volume
 one in $\BBR^d$.

\begin{itemize}

\item \hbox{\rm ( Expectation)} There is a positive constant $c_K$
such that  $$\E(\Vol(K_n)) = 1- (c_K+o(1)) n^{-2/(d+1)}. $$

\item \hbox{\rm (Variance)} $\Var (\Vol(K_n)) = \Theta (n^{-(d+3)/(d+1)}).
$

\item \hbox{\rm (Higher Moments)} For any fixed  $k$, $M_k (\Vol
(K_n)) = \Theta (n^{-k +(d-1)/(d+1)})$. For any fixed odd $k$,
absolute value of the $k$th moment of $\Vol (K_n)$ is $$O( n^{-k
+(d-2 +o(1))/(d+1)}) = O(M_k n^{-1/(d+1)+o(1)}). $$

\item \hbox{\rm (Speed of Convergence)} Almost surely,

$$ \lim_{n \rightarrow \infty} |\frac{\Vol(K_n)_n} {\BE(\Vol(K_n)_n)} -1| \delta (n)
n^{(d+3)/(d+1)}  \ln^{-1/2} n =0 $$

\noindent for any function $\delta (n)$ tending to zero with $n$.

\item \hbox{\rm (Central Limit Theorem)} $\Vol(K_n)$ satisfies the central
limit theorem.

\item \hbox{\rm (Exponential Tail)} There is a positive constant
$c$ such that for any $$0 < \lambda \le \alpha
n^{\frac{(d-1)(d+3)}{(d+1)(3d+5)}}$$ we have

$$  \BP(|\Vol(K_n)-\BE(\Vol(K_n))| \ge \sqrt {\lambda
\Var(\Vol(K_n)) }) \le 2\exp(-c\lambda) + \exp (-c
n^{(d-1)/(3d+5)}). $$

\end{itemize} \end{theorem}

\begin{remark} The  statement on the expectation is due to
B\'arany \cite{Bar3} and Sch\"utt \cite{Sch}. The  statement on
the variance is due to  Reitzner \cite{Rei, Rei2}. For the upper
bound on the variance, Vu \cite{Vu1} has a different proof, which
also extends to  $M_k$ for any fixed $k$. The lower bound for
$M_k$ follows from the central limit theorem. The estimate for odd
moments follows the bound on $M_k$ and the error term $\ep_n$ in
Theorem \ref{theo:CLT}. The statement concerning the central limit
theorem is Theorem \ref{theo:CLT}. The statement concerning
exponential tail follows from Theorem \ref{theo:Vu1}. The
statement on the speed of convergence is a corollary of this
theorem.
\end{remark}

 \vskip2mm

\noindent {\it Facets.}

\begin{theorem} \label{theo:sumarizefacet} Let $K$ be a smooth convex body of volume
 one in $\BBR^d$ and $0 \le i < d$ be an integer.

\begin{itemize}

\item \hbox{\rm ( Expectation)} There is a positive constant $c_i$
such that  $$\E(f_i) =  (c_i+o(1)) n^{(d-1)/(d+1)}. $$

\item \hbox{\rm (Variance)} $\Var (f_i) = \Theta (n^{(d-1)/(d+1)}).
$

\item \hbox{\rm (Higher Moments)} For any fixed $k$, $M_k (f_i) =
\Theta (n^{k(d-1)/2(d+1)}). $ For any fixed odd $k$, the absolute
value of the $k$th moment of $f_i$ is $$O(
n^{(k(d-1)-2)/2(d+1)+o(1) }) = O(M_k n^{-1/(d+1)+o(1)}). $$

\item \hbox{\rm (Speed of Convergence)} Almost surely,

$$ \lim_{n \rightarrow \infty} |\frac{f_i} {\BE(f_i)} -1| \delta (n)
n^{-(d-1)/(d+1)}  \ln^{-1/2} n =0 $$

\noindent for any function $\delta (n)$ tending to zero with $n$.

\item \hbox{\rm (Central Limit Theorem)} $f_i$ satisfies the central
limit theorem.

\item \hbox{\rm (Exponential Tail)} There are positive constants $c_i, \gamma_i, \beta_i$  such that for any $$0 < \lambda \le
\gamma_i n^{\beta_i}$$ we have

$$ \BP(|f_i-\BE(f_i)| \ge \sqrt {\lambda \Var(f_i) }) \le
2\exp(-c_i\lambda). $$

\end{itemize} \end{theorem}

\begin{remark} The asymptotic of the expectation is from
Reitzner \cite{Rei3}. Two special cases $i=0$ and $i=d-1$ were
established earlier by B\'ar\'any \cite{Bar3} and  Weieacker
\cite{Wei}, respectively. The origin of the remaining statements
is the same as in the previous remark.
\end{remark}

\section{More central limit theorems}

As the reader has possibly  noticed, the  key ingredients in our
approach are: the  results in Section 2 concerning the Poisson
model, the concentration results and the bounds on the
differences. Once these ingredients are available, we are ready to
prove a central limit theorem.

A nice thing about the current method is that the proofs involving
in these ingredients are quite robust, at least  in spirit. They
involve very abstract tools which have little to do with the
specific model of the random polytope one deals with.
Representative examples are the Baldi-Rinott theorem used in
\cite{Rei} and the Divide and Conquer Martingale technique used in
\cite{Vu1}. Thus, typically we can extend the results obtained
here  for many other models of random polytopes (the actual
argument may get complicated in some occasions). For instance,
analogues of the results in Section 2 are now available for the
case when the body $K$ is a polytope, due to a recent result of
B\'ar\'any and Reitzner (private communication from B\'ar\'any).
Furthermore, concentration results are ready for this model as
well (see \cite{Vu1}). Thus one  can prove the corresponding
central limit theorems. Details will appear elsewhere.

\section {Appendix}

{\it \noindent Proof of the second half of Lemma
\ref{lemma:main3}.} Consider a floating body $F= F_{c_0 \log n/n}$
for some constant $c_0$. An observation from \cite{Vu1} shows that
there is a collection of $m= O(n^{c_1})$ caps $C_1, \dots, C_m$,
each of volume $c_2 \log n/n$, where both $c_1$ and $c_2$ are
constants, such that the following holds: Let $x$ be a point in $K
\backslash F$. Then the region seen by $x$ is contained in one of
the $C_i$.

Now let us fix a constant $C$ as in Lemma \ref{lemma:main3}.
Choose $c_0$ sufficiently large  so that the probability that
$F_{c_0 \log n/n} \subset K_n$ is at least $1 - n^{-C-1}$. The
choice of $c_0$ determines $c_1$ and $c_2$. Next, we choose $c_5$
sufficiently large so that the probability that each of the $C_i$
contains at least $c_5 \log n$ points is at most $n^{-C-c_1-1})$.

We say that a set $P$ of $n$ random points is regular if $F_{c_0
\log n/n} \subset \Conv (P)$ and each $C_i$ contains less than
$c_5 \log n$ points from $P$.

It is clear that the probability that there is no point which is
$c_5 \log n$-wide with respect to $P$ is at least  the probability
that $P$ is regular, which is bounded from below by

$$ 1-n^{-C-1} -  m n^{-C-c_1-1} =  1-n^{-C-1} -  n^{c_1}
n^{-C-c_1-1} \ge 1- n^{-C}, $$

\noindent completing the proof. \hs

\end{document}